\documentclass[12pt]{article}

\usepackage{amsmath, epsfig, cite}
\usepackage{amssymb}
\usepackage{amsfonts}
\usepackage{latexsym}

\newtheorem{thm}{Theorem}[section]

\newtheorem{conj}[thm]{Conjecture}
\newtheorem{lem}[thm]{Lemma}

\newcommand{\pf}{\noindent{\it Proof.} }
\newcommand*{\fl}[2]{\left\lfloor\frac{#1}{#2}\right\rfloor}

\def\fl#1{\left\lfloor#1\right\rfloor}

\numberwithin{equation}{section}

\newcommand{\qed}{{\hfill$\square$}\medskip}

\begin{document}

\nocite{*}
\begin{center}
{\Large\bf Proof of a conjecture involving Sun polynomials}
\end{center}

\vskip 2mm \centerline{Victor J. W. Guo$^1$, Guo-Shuai Mao$^2$, Hao Pan$^3$}
\begin{center}
{\footnotesize $^1$School of Mathematical Sciences, Huaiyin Normal University, Huai'an, Jiangsu 223300,
 People's Republic of China\\
{\tt jwguo@hytc.edu.cn}

\vskip 2mm

$^{2,3}$Department of Mathematics, Nanjing University, Nanjing 210093, People's Republic of China\\
{\tt mg1421007@smail.nju.edu.cn, haopan79@zoho.com} }

\end{center}


\vskip 0.7cm \noindent{\bf Abstract.} The Sun polynomials $g_n(x)$ are defined by
\begin{align*}
g_n(x)=\sum_{k=0}^n{n\choose k}^2{2k\choose k}x^k.
\end{align*}
We prove that, for any positive integer $n$, there hold
\begin{align*}
&\frac{1}{n}\sum_{k=0}^{n-1}(4k+3)g_k(x) \in\mathbb{Z}[x],\quad\text{and}\\
&\sum_{k=0}^{n-1}(8k^2+12k+5)g_k(-1)\equiv 0\pmod{n}.
\end{align*}
The first one confirms a recent conjecture of Z.-W. Sun, while the second one partially answers another conjecture of Z.-W. Sun.
We give three different proofs of the former. One of them depends on the following congruence:
$$
{m+n-2\choose m-1}{n\choose m}{2n\choose n}\equiv 0\pmod{m+n}\quad\text{for $m,n\geqslant 1$.}
$$

\vskip 3mm \noindent {\it Keywords}: congruence; Sun polynomials; Chu-Vandermonde's identity; $q$-binomial coefficients; reciprocal and unimodal polynomials;
Zeilberger's algorithm.

\vskip 2mm
\noindent{\it MR Subject Classifications}: 11A07, 11B65, 05A10

\section{Introduction}
Recently, Z.-W. Sun \cite{Sun} introduced the polynomials
\begin{align*}
g_n(x)=\sum_{k=0}^n{n\choose k}^2{2k\choose k}x^k,
\end{align*}
which we call {\it Sun polynomials} here, and proved many interesting identities and congruences involving $g_n(x)$, such as
\begin{align*}
\frac{1}{3n^2}\sum_{k=0}^{n-1}(4k+3)g_k
&=\sum_{k=0}^{n-1}\frac{1}{k+1}{2k\choose k}{n-1\choose k}^2,\\[5pt]
\sum_{k=0}^{p-1}kg_k
&\equiv -\frac{3}{4} \pmod{p^2},
\end{align*}
where $p$ is an odd prime and $g_k=g_k(1)$. Z.-W. Sun \cite{Sunnew} also conjectured that
\begin{align*}
\sum_{k=0}^\infty \frac{16k+5}{324^k}{2k\choose k}g_k(-20)=\frac{189}{25\pi}.
\end{align*}
 Some other congruences involving $g_n$ can be found in \cite{MS,Sun13,Sunnew}.

The Sun polynomials also satisfy the following identities\cite[(2.7),(2.11)]{Sun}:
\begin{align*}
\sum_{k=0}^n {n\choose k}f_k(x) &=g_n(x),\\[5pt]
\sum_{k=0}^n {n\choose k}{n+k\choose k}(-1)^{n-k}g_k(x) &=A_n(x),
\end{align*}
where $f_n(x)$ and $A_n(x)$ are respectively the Franel polynomials and Ap\'ery polynomials \cite{Sun00} defined as
\begin{align*}
f_n(x)&=\sum_{k=0}^{n}{n\choose k}^2{2k\choose n}x^k,\\[5pt]
A_n(x)&=\sum_{k=0}^{n}{n\choose k}^2 {n+k\choose k}^2 x^k.
\end{align*}

The objective of this paper is to prove the following result, which was originally conjectured by Z.-W. Sun (see \cite[Conjecture 4.1(ii)]{Sun}).
\begin{thm}\label{thm:1}
Let $n$ be a positive integer. Then
\begin{align}
&\frac{1}{n}\sum_{k=0}^{n-1}(4k+3)g_k(x) \in\mathbb{Z}[x], \label{eq:first}\\
&\sum_{k=0}^{n-1}(8k^2+12k+5)g_k(-1)\equiv 0\pmod{n}, \label{eq:second}.
\end{align}
\end{thm}

\noindent{\it Remark.} For the congruence \eqref{eq:second}, Z.-W. Sun \cite{Sun} made the following stronger conjecture:
\begin{align}
\sum_{k=0}^{n-1}(8k^2+12k+5)g_k(-1) &\equiv n^2\pmod{2n^2}, \label{eq:second-2}\\
\sum_{k=0}^{p-1}(8k^2+12k+5)g_k(-1) &\equiv 3p^2\pmod{p^3}, \label{eq:third-2}
\end{align}
where $p$ is a prime.

In order to prove Theorem \ref{thm:1}, we need to establish some preliminary results in Section~\ref{sec:lemmas}.
However, since the following result is interesting in its own right, we label it as a theorem here.

\begin{thm}\label{thm:2}
Let $m$ and $n$ be positive integers. Then
\begin{align}
{m+n-2\choose m-1}{n\choose m}{2n\choose n}\equiv 0\pmod{m+n}. \label{eq:important}
\end{align}
\end{thm}

It is worth mentioning that Gessel \cite[Section 7]{Gessel} proved a similar result as follows:
$$
\frac{m}{2}{2m\choose m}{2n\choose n}\equiv 0\pmod{m+n},
$$
of which a generalization was given by the author \cite[Theorem~1.4]{Guo}.

The paper is organized as follows. Applying the same techniques in
\cite{Guo,GK}, we shall prove a $q$-analogue of Theorem \ref{thm:2} in the next section. In Section~3, we give three lemmas, one of which
is closely related to Theorem \ref{thm:2}. Two proofs of \eqref{eq:first} and a proof of \eqref{eq:second} will be given in Section 4.
The second proof of \eqref{eq:first} is motivated by Sun \cite[Lemma 3.4]{Sun} and its proof. We shall also give a $q$-analogue (the third proof)
of \eqref{eq:first} in Section 5. We end the paper in Section 6 with a related conjecture.

\section{A $q$-analogue of Theorem \ref{thm:2}}
Recall that the {\it $q$-binomial coefficients}
are defined by
\begin{align*}
{n\brack k}_q
=\begin{cases}
\displaystyle \prod_{i=1}^{k}\frac{1-q^{n-k+i}}{1-q^i},
&\text{if $0\leqslant k\leqslant n,$} \\[10pt]
0, &\text{otherwise.}
\end{cases}
\end{align*}
We now state the announced strengthening of Theorem \ref{thm:2}.

\begin{thm}\label{thm:q-analog}
Let $m$ and $n$ be positive integers. Then
\begin{align}
\frac{1-q}{1-q^{m+n}}{m+n-2\brack m-1}_q{n\brack m}_q{2n\brack n}_q \label{eq:q-analog}
\end{align}
is a polynomial in $q$ with non-negative integer coefficients.
\end{thm}

It is clear that Theorem~\ref{thm:2}
can be deduced from Theorem~\ref{thm:q-analog} by letting $q\to 1$.

A polynomial
$A(q)=
\sum _{i=0} ^{d}a_iq^i$ in $q$ of degree $d$ is called {\it reciprocal\/} if
$a_i=a_{d-i}$ for all $i$, and that it
is called {\it unimodal\/} if there is an index $r$ such that $0\leqslant a_0\leqslant\dots
\leqslant a_r\geqslant\dots\geqslant a_d\geqslant 0$.
The following is an elementary but crucial property of reciprocal and unimodal polynomials
(see, for example, \cite{Andrews75} or \cite[Proposition 1]{Stanley89}).

\begin{lem}\label{lem:AqBq}
If $A(q)$ and $B(q)$ are reciprocal and unimodal polynomials, then so is
their product $A(q)B(q)$.
\end{lem}

Similarly to the proof of \cite[Theorem 3.1]{GK}, we also need the following result.
We refer the reader to \cite[Proposition~10.1.(iii)]{ReSWAA} and \cite[Proof of Theorem~2]{AndrCB}
for similar mathematical ideas.

\begin{lem}{\rm\cite[Lemma 5.1]{GK}} \label{lem:RSW}
Let $P(q)$ be a reciprocal and unimodal
polynomial and $m$ and $n$ positive integers with $m\leqslant n$.
Furthermore, assume that $\frac {1-q^m} {1-q^n}P(q)$ is a polynomial
in $q$. Then $\frac {1-q^m} {1-q^n}P(q)$ has non-negative coefficients.
\end{lem}

\noindent{\it Proof of Theorem {\rm\ref{thm:q-analog}.}}
It is well known that the $q$-binomial coefficients are reciprocal and unimodal
polynomials in $q$ (see, for example, \cite[Ex.~7.75.d]{Stanley}). By Lemma \ref{lem:AqBq}, so is
the product of three $q$-binomial coefficients. In view of Lemma~\ref{lem:RSW}, to prove Theorem~\ref{thm:q-analog}, it
suffices to show that the expression \eqref{eq:q-analog} is a polynomial in $q$. We shall accomplish this by considering
a count of cyclotomic polynomials.

Recall that
$$
q^n-1=\prod _{d\mid n} ^{}\Phi_d(q),
$$
where $\Phi_d(q)$ denotes the $d$-th cyclotomic polynomial in $q$.
Therefore,
$$
\frac{1-q}{1-q^{m+n}}{m+n-2\brack m-1}_q{n\brack m}_q{2n\brack n}_q
=\prod _{d=2} ^{2n}\Phi_d(q)^{e_d},
$$
with
\begin{align*}
e_d&
={}-\chi(d\mid m+n)
+\fl{\frac {m+n-2} {d}}
+\fl{\frac {2n} {d}} \notag\\
&\quad{}-\fl{\frac {m-1} {d}}
-\fl{\frac {n-1} {d}}
-\fl{\frac {m} {d}}-
\fl{\frac {n} {d}}
-\fl{\frac {n-m} {d}}, 
\end{align*}
where $\chi(\mathcal S)=1$ if $\mathcal S$ is
true and $\chi(\mathcal S)=0$ otherwise.
The number $e_d$ is obviously non-negative, unless
$d\mid m+n$.

So, let us assume that $d\mid m+n$ and $d\geqslant 2$. We consider two cases: If $d\mid m$, then $d\mid n$, and so
\begin{align*}
\left\lfloor\frac{m+n-2}{d}\right\rfloor-\left\lfloor\frac{m-1}{d}\right\rfloor-\left\lfloor\frac{n-1}{d}\right\rfloor
=\frac{m+n-d}{d}-\frac{m-d}{d}-\frac{n-d}{d}=1.
\end{align*}
Namely, $e_d=0$ is non-negative; If $d\nmid m$, then
\begin{align*}
\left\lfloor\frac{2n}{d}\right\rfloor-\left\lfloor\frac{m}{d}\right\rfloor-\left\lfloor\frac{n}{d}\right\rfloor
-\left\lfloor\frac{n-m}{d}\right\rfloor
=\frac{m+n}{d}-\left\lfloor\frac{m}{d}\right\rfloor-\left\lfloor\frac{n}{d}\right\rfloor
=1.
\end{align*}
That is, $e_d=0$ is still non-negative.
This completes the proof of polynomiality of \eqref{eq:q-analog}.
\qed

\section{Some preliminary results} \label{sec:lemmas}

\begin{lem}\label{lem:one}
Let $n$ be a non-negative integer. Then
\begin{align}
{x\choose n}^2=\sum_{k=0}^{n}{x\choose n+k}{n+k\choose k}{n\choose k}. \label{eq:chu-van}
\end{align}
\end{lem}
\pf  Applying Chu-Vandermonde's identity (see, for example, \cite[p.~32]{Koepf})
\begin{align*}
{x\choose n}=\sum_{k=0}^{n}{x-n\choose k}{n\choose k},
\end{align*}
and noticing that ${x\choose n}{x-n\choose k}={x\choose n+k}{n+k\choose k}$, we obtain \eqref{eq:chu-van}.
In fact, Eq.~\eqref{eq:chu-van} is a special case of \cite[p.~15, Eq.~(9)]{Riordan}.
\qed

\begin{lem}\label{lem:two}
Let $n$ be a positive integer and let $0\leqslant k\leqslant n$. Then
\begin{align*}
\sum_{m=k}^{n-1}(4m+3){m\choose k}&=(4n-1){n\choose k+1}-4{n\choose k+2}, \\[5pt]
\sum_{m=k}^{n-1}(8m^2+12m+5){m\choose k}&=(8n^2-4n+1){n\choose k+1}-(16n-12){n\choose k+2}+16{n\choose k+3}.
\end{align*}
\end{lem}
\pf Proceed by induction on $n$.
\qed

\begin{lem}\label{lem:three}
Let $m$ and $n$ be non-negative integers. Then
\begin{align*}
{m+n\choose m}{n+1\choose m}{2n\choose n}\frac{3m^2+n^2+m+n}{(m+n)(n+1)}\equiv 0\pmod{m+n+1}.
\end{align*}
\end{lem}
\pf Note that
\begin{align*}
\frac{3m^2+n^2+m+n}{(m+n)(n+1)}=\frac{3m(m+n+1)}{(m+n)(n+1)}+\frac{n(m+n+1)}{(m+n)(n+1)}-\frac{2m(2n+1)}{(m+n)(n+1)}.
\end{align*}
Since $\frac{1}{n+1}{2n\choose n}={2n\choose n}-{2n\choose n-1}$ is an integer (the $n$-th Catalan number), we see that
\begin{align*}
{m+n\choose m}{2n\choose n}\frac{m}{(m+n)(n+1)}={m+n-1\choose m-1}{2n\choose n}\frac{1}{n+1},
\end{align*}
and
\begin{align*}
{m+n\choose m}{2n\choose n}\frac{n}{(m+n)(n+1)}={m+n-1\choose m}{2n\choose n}\frac{1}{n+1},
\end{align*}
are both integers. It remains to show that
\begin{align*}
{m+n\choose m}{n+1\choose m}{2n\choose n}\frac{2m(2n+1)}{(m+n)(n+1)}\equiv 0\pmod{m+n+1},
\end{align*}
i.e.,
\begin{align*}
{m+n-1\choose m-1}{n+1\choose m}{2n+2\choose n+1}\equiv 0\pmod{m+n+1}.
\end{align*}
But this is just the $n\to n+1$ case of the congruence \eqref{eq:important}.
\qed

\begin{lem}
Let
\begin{align*}
S_n=f_{n-3}(-1)=\sum_{k=0}^{n-3}(-1)^k{2k\choose k}{n-3\choose k}{k\choose n-k-3}. 
\end{align*}
Then there hold the following congruences:
\begin{align}
S_{3n}&\equiv S_{3n+1}\equiv -S_{3n+2}\pmod{3}, \label{eq:rec-10}\\
S_{4n+2}&\equiv 0\pmod{4}, \label{eq:rec-11}\\
S_{n+2}+12S_{n+1}+16S_n &\equiv 0\pmod{n}. \label{eq:rec-1}
\end{align}
\end{lem}
\pf
Zeilberger's algorithm \cite{Koepf,PWZ} gives the following recurrence relation for $S_n$:
\begin{align}
&(5n^3-8n^2)S_{n+3}+(45n^3-117n^2+90n-24)S_{n+2} \notag\\
&\quad{}+(200n^3-720n^2+824n-288)S_{n+1}+(160n^3-736n^2+1024n-384)S_n=0.
\label{eq:rec}
\end{align}
Replacing $n$ by $3n-1$ in \eqref{eq:rec}, we obtain
$$-S_{3n+2}-S_{3n}\equiv 0\pmod 3,$$
while replacing $n$ by $3n+1$ in \eqref{eq:rec}, we get
$$S_{3n+2}+S_{3n+1}\equiv 0\pmod 3.$$
This proves \eqref{eq:rec-10}. Similarly, replacing $n$ by $4n-1$ in \eqref{eq:rec}, we are led to \eqref{eq:rec-11}.

In order to prove \eqref{eq:rec-1}, we need to consider four cases:
\begin{itemize}
\item If $\gcd(n,24)=1$, then \eqref{eq:rec} means that $-24(S_{n+2}+12S_{n+1}+16S_n)\equiv 0\pmod{n}$, i.e.,
the congruence \eqref{eq:rec-1} holds.
\item If $\gcd(n,24)=2,4,8$, then \eqref{eq:rec} means that
$$\pm 2nS_{n+2}-24(S_{n+2}+12S_{n+1}+16S_n)\equiv 0\pmod{8n}.$$
By \eqref{eq:rec-11}, we have $2nS_{n+2}\equiv 0\pmod{8n}$ in this case, and so the congruence \eqref{eq:rec-1} holds.
\item If $\gcd(n,24)=3$, then \eqref{eq:rec} means that
$$2nS_{n+1}+nS_n-24(S_{n+2}+12S_{n+1}+16S_n)\equiv 0\pmod{3n}.$$
By \eqref{eq:rec-10}, we have $2nS_{n+1}+nS_n\equiv 0\pmod{3n}$ in this case, and so
the congruence \eqref{eq:rec-1} holds.
\item If $\gcd(n,24)=6,12,24$, then \eqref{eq:rec} means that
$$30nS_{n+2}+8nS_{n+1}+16nS_n-24(S_{n+2}+12S_{n+1}+16S_n)\equiv 0\pmod{24n},\quad\text{or}$$
$$18nS_{n+2}+8nS_{n+1}+16nS_n-24(S_{n+2}+12S_{n+1}+16S_n)\equiv 0\pmod{24n}.$$
By \eqref{eq:rec-11}, we have $30nS_{n+2}\equiv 18nS_{n+2}\equiv 0\pmod{24n}$ and $8nS_{n+1}+16nS_n\equiv 0\pmod{24n}$ in this case,
and so the congruence \eqref{eq:rec-1} still holds.
\end{itemize}

\qed

\section{Proof of Theorem \ref{thm:1} }
\noindent{\it First Proof of \eqref{eq:first}.}
By Lemmas \ref{lem:one} and \ref{lem:two}, we have
\begin{align}
&\hskip -2mm \sum_{m=0}^{n-1}(4m+3)g_m(x) \notag\\[5pt]
&=\sum_{m=0}^{n-1}(4m+3)\sum_{k=0}^{m}{m\choose k}^2{2k\choose k}x^k  \notag\\[5pt]
&=\sum_{m=0}^{n-1}(4m+3)\sum_{k=0}^{m}{2k\choose k}x^k\sum_{i=0}^k{m\choose k+i}{k+i\choose i}{k\choose i} \notag\\[5pt]
&=\sum_{k=0}^{n-1}{2k\choose k}x^k\sum_{i=0}^k {k+i\choose i}{k\choose i} \sum_{m=k+i}^{n-1}(4m+3){m\choose k+i} \notag\\[5pt]
&=\sum_{k=0}^{n-1}{2k\choose k}x^k\sum_{i=0}^k \left((4n-1){n\choose k+i+1}-4{n\choose k+i+2}\right){k+i\choose i}{k\choose i}. \label{eq:multi-sum}
\end{align}
For any non-negative integer $k\leqslant n-1$, to prove that the coefficient of $x^k$ in the right-hand side of \eqref{eq:multi-sum}
is a multiple of $n$, it suffices to show that
\begin{align}
{2k\choose k}\sum_{i=0}^k \left({n\choose k+i+1}+4{n\choose k+i+2}\right){k+i\choose i}{k\choose i}\equiv 0\pmod{n}. \label{eq:single-sum}
\end{align}

We shall accomplish the proof of \eqref{eq:single-sum} by using a minor trick.  Rewrite the left-hand side of \eqref{eq:single-sum} as
\begin{align*}
&\hskip -2mm {2k\choose k}\sum_{i=0}^{k+1} {n\choose k+i+1}\left({k+i\choose i}{k\choose i}+4{k+i-1\choose i-1}{k\choose i-1}\right) \\
&= \sum_{i=0}^{k+1} {n\choose k+i+1}{k+i\choose i}{k+1\choose i}{2k\choose k}\frac{k^2+3i^2+k+i}{(k+i)(k+1)}.
\end{align*}
Then, by Lemma \ref{lem:three}, for each $i\leqslant k+1$, the expression
$$
{k+i\choose i}{k+1\choose i}{2k\choose k}\frac{k^2+3i^2+k+i}{(k+i)(k+1)}.
$$
is a multiple of $k+i+1$. Finally, noticing that $${n\choose k+i+1}(k+i+1)=n{n-1\choose k+i}\equiv 0\pmod{n},$$
we complete the proof.
\qed

\noindent{\it Second Proof of \eqref{eq:first}.} This proof is motivated by \cite[Lemma 3.4 and its proof]{Sun}.
It is clear that \eqref{eq:first} is equivalent to the following congruence:
\begin{align}
{2j\choose j}\sum_{k=j}^{n-1}(4k+3){k\choose j}^2\equiv 0\pmod n. \label{eq:mao}
\end{align}
Denote the left-hand side of \eqref{eq:mao} by $u_j$. Then by Zeilberger's algorithm \cite{PWZ}, we have
\begin{align}
u_{j+1}-u_j&=-{2j\choose j}{n-1\choose j}^2\frac{(9j+6)(j+1)n^2+(12j^2-8jn-4n+14j+4)n^3}{(j+1)^3(j+2)} \notag\\[5pt]
&=-{2j\choose j}{n-1\choose j}{n+1\choose j+1}\frac{(9j+6)n^2}{(j+1)(n+1)} \notag\\[5pt]
&\quad{}-{2j\choose j}{n\choose j+1}{n+1\choose j+1}\frac{(12j^2-8jn-4n+14j+4)n^2}{(j+1)(n+1)}. \label{eq:uj-rec}
\end{align}
Noticing that $\frac{1}{j+1}{2j\choose j}$ is an integer and $n+1$ is relatively prime to $n$, from \eqref{eq:uj-rec} we immediately get
$$ u_{j+1}-u_j\equiv0\pmod n.$$
Since $u_0=2n^2+n\equiv0\pmod n$, we conclude that $u_j\equiv0\pmod n$ for all $j$.
This proves \eqref{eq:mao}.

\medskip
\noindent{\it Proof of \eqref{eq:second}.}
By Lemmas \ref{lem:one} and \ref{lem:two}, similarly to \eqref{eq:multi-sum}, we have
\begin{align}
&\hskip -2mm \sum_{m=0}^{n-1}(8m^2+12m+5)g_m(-1) \notag\\[5pt]
&=\sum_{k=0}^{n-1}{2k\choose k}(-1)^k\sum_{i=0}^k \Bigg((8n^2-4n+1){n\choose k+i+1}-(16n-12){n\choose k+i+2} \notag\\
&\qquad{}+16{n\choose k+i+3}\Bigg){k+i\choose i}{k\choose i}. \label{eq:multi-sum-2}
\end{align}
In view of \eqref{eq:single-sum}, it follows from \eqref{eq:multi-sum-2} that
\begin{align}
&\hskip -2mm
\sum_{m=0}^{n-1}(8m^2+12m+5)g_m(-1) \notag\\[5pt]
&\equiv \sum_{k=0}^{n-1}{2k\choose k}(-1)^k\sum_{i=0}^k \left({n\choose k+i+1}+12{n\choose k+i+2}+16{n\choose k+i+3}\right) \notag \\
&\qquad{}\times{k+i\choose i}{k\choose i} \notag \\
&=\sum_{m=0}^{n-1} \left({n\choose m+1}+12{n\choose m+2}+16{n\choose m+3}\right) \notag\\
&\qquad{}\times\sum_{k=0}^{m}(-1)^k{2k\choose k}{m\choose k}{k\choose m-k} \pmod{2n^2}. \label{eq:multi-sum-3}
\end{align}
Note that the right-hand side of \eqref{eq:multi-sum-3} may be written as
\begin{align}
&\hskip -2mm\sum_{m=0}^{n-1} \left({n\choose m+1}+12{n\choose m+2}+16{n\choose m+3}\right) S_{m+3} \notag\\
&=\sum_{m=1}^{n}{n\choose m}(S_{m+2}+12S_{m+1}+16S_m), \label{eq:rewrite}
\end{align}
which is clearly congruent to $0$ modulo $n$ by \eqref{eq:rec-1} and the fact that $m{n\choose m}=n{n-1\choose m-1}$.  \qed

\section{A $q$-analogue of \eqref{eq:first}}
Define the $q$-analogue of Sun polynomials as follows:
$$
g_n(x;q)=\sum_{k=0}^n {n\brack k}_{q}^2 {2k\brack k}_{q}x^k.
$$
We have the following congruences related to $g_n(x;q)$.
\begin{thm}\label{thm:5}
Let $n$ be a positive integer. Then
\begin{align}
&(1+q)^2\sum_{k=0}^{n-1}q^{2k}[k+1]_{q^{2}}g_k(x;q^2)
\equiv \sum_{k=0}^{n-1}q^{2k}g_k(x;q^2)\pmod{\prod_{\substack{d\mid n\\ d>1\text{ is odd}}}\Phi_d(q)}, \label{oddcong}  \\
&\sum_{k=0}^{n-1}q^{k}g_k(x;q)
\equiv 0\pmod{\prod_{\substack{d\mid n\\ d\text{ is even}}}\Phi_d(q)}, \label{evencong1}\\
&\sum_{j=0}^{n-1}x^j [j+1]_{q^2}{2j\brack j}_q\sum_{k=j}^{n-1}q^k{k\brack j}_q {k+1\brack j+1}_q
\equiv 0\pmod{\prod_{\substack{d\mid n\\ d>2\text{ is even}}}\Phi_d(q)}, \label{evencong2}
\end{align}
where $[n]_q=\frac{1-q^n}{1-q}$ denotes a $q$-integer. 
\end{thm}

\pf It is clear that
\begin{align*}
\sum_{k=0}^{n-1}q^{2k}[k+1]_{q^{2}}g_k(x;q^2)&=
\sum_{k=0}^{n-1}q^{2k}[k+1]_{q^{2}}\sum_{j=0}^k {k\brack j}_{q^2}^2{2j\brack j}_{q^2}x^j \\
&=\sum_{j=0}^{n-1}[j+1]_{q^2}{2j\brack j}_{q^2}x^j\sum_{k=j}^{n-1}q^{2k}
{k+1\brack j+1}_{q^2}{k\brack j}_{q^2}.
\end{align*}
Suppose that $d\mid n$ and $d$ is odd. It is easy to see that $\Phi_d(q)$ divides $\Phi_d(q^2)$.
Write $j=\gamma d+\delta$, where $0\leqslant \delta\leqslant d-1$. If $d\leqslant 2\delta$, then by the $q$-Lucas theorem (see Olive \cite{Olive},
D\'esarm\'enien \cite[Proposition 2.2]{Des} or Guo and Zeng \cite[Proposition 2.1]{GZ}),
$$
{2j\brack j}_{q^2}\equiv {2\gamma+1\choose \gamma}{2\delta-d\brack \delta}_{q^2}=0\pmod{\Phi_d(q)}.
$$
Now assume that $\delta\leqslant \frac{d-1}{2}$.
Then applying the $q$-Lucas theorem, we have
\begin{align}\label{oddqlucas}
\sum_{k=j}^{n-1}q^{2k}
{k+1\brack j+1}_{q^2} {k\brack j}_{q^2}=&\sum_{\alpha=0}^{\frac{n}{d}-1}\sum_{\beta=0}^{d-1}
q^{2(\alpha d+\beta)}
{\alpha d+\beta+1\brack \gamma d+\delta+1}_{q^2} {\alpha d+\beta\brack \gamma d+\delta}_{q^2}\notag\\
\equiv&\sum_{\alpha=0}^{\frac{n}{d}-1}\binom{\alpha}{\gamma}^2
\sum_{\beta=0}^{d-1}q^{2\beta} {\beta+1\brack \delta+1}_{q^2} {\beta\brack \delta}_{q^2}\pmod{\Phi_d(q)}.
\end{align}

It is easy to see that
\begin{align}\label{oddcv}
\sum_{\beta=0}^{d-1}
q^{2\beta}
{\beta+1 \brack \delta+1}_{q^2} {\beta\brack \delta}_{q^2}
&=\sum_{r=0}^{d-1-\delta}q^{2r+2\delta}
{r+\delta+1\brack \delta+1}_{q^2} {r+\delta\brack \delta}_{q^2}\notag\\
&=\sum_{r=0}^{d-1-\delta}q^{2\delta+6r+4r\delta+2r^2}
{-\delta-2\brack r}_{q^2} {-\delta-1\brack r}_{q^2}\notag\\
&\equiv \sum_{r=0}^{d-1-\delta}q^{2(d-2-\delta-r)(d-1-\delta-r)-4\delta-2\delta^2-4}
{d-\delta-2 \brack r}_{q^2} {d-\delta-1\brack d-\delta-1-r}_{q^2}\notag\\
&= q^{-4\delta-2\delta^2-4}{2d-2\delta-3\brack d-1-\delta}_{q^2}\pmod{\Phi_{d}(q)},
\end{align}
where we have used the $q$-Chu-Vandemonde identity (see \cite[(3.3.10)]{Andrews98}) in the last step. 
Furthermore, we have
$$
{2d-2\delta-3\brack d-1-\delta}_{q^2}\equiv 0\pmod{\Phi_d(q)}
$$
for $\delta\leqslant \frac{d-3}{2}$.  This proves that
\begin{align*}
&\hskip -2mm
\sum_{k=j}^{n-1}q^{2k}
{k+1\brack j+1}_{q^2} {k\brack j}_{q^2} \\
&\equiv \begin{cases} \displaystyle \sum_{\alpha=0}^{\frac{n}{d}-1}\binom{\alpha}{\gamma}^2 q^{\frac{d-5}{2}}{d-2\brack \frac{d-1}{2}}_{q^2},&\text{if $j\equiv \frac{d-1}{2}\pmod d$,}
\\[10pt]
0,&\text{otherwise.}
\end{cases}\pmod{\Phi_d(q)}.
\end{align*}
Hence, writing $j=\gamma d+\frac{d-1}{2}$ and applying the $q$-Lucas theorem, we obtain
\begin{align*}
&\hskip -2mm
\sum_{j=0}^{n-1}[j+1]_{q^2} {2j\brack j}_{q^2}x^j\sum_{k=j}^{n-1}q^{2k}
{k+1\brack j+1}_{q^2} {k\brack j}_{q^2}\\
&\equiv \sum_{\gamma=0}^{\frac{n}{d}-1}\sum_{\alpha=0}^{\frac{n}{d}-1}{2\gamma\choose\gamma}{\alpha\choose\gamma}^2x^{\gamma d+\frac{d-1}{2}}
q^{\frac{d-5}{2}}\left[\frac{d+1}{2}\right]_{q^2} {d-1\brack \frac{d-1}{2}}_{q^2} {d-2\brack \frac{d-1}{2}}_{q^2}\\
&\equiv \sum_{\gamma=0}^{\frac{n}{d}-1}\sum_{\alpha=0}^{\frac{n}{d}-1}{\alpha\choose\gamma}^2
\frac{x^{\gamma d+\frac{d-1}{2}}}{q(1+q)^2}\pmod{\Phi_d(q)},
\end{align*}
where we have used the congruence
$$
{d-1\brack k}_{q^2}\equiv(-1)^k q^{-k(k+1)}\pmod{\Phi_d(q)}\quad\text{for }0\leqslant k\leqslant d-1.
$$

On the other hand, we have
\begin{align*}
\sum_{k=0}^{n-1}q^{2k}g_k(x;q^2) 
=\sum_{j=0}^{n-1} {2j\brack j}_{q^2}x^j\sum_{k=j}^{n-1}q^{2k} {k\brack j}_{q^2}^2.
\end{align*}
Similarly as before, if $j=\gamma d+\delta$ and $\delta\leqslant \frac{d-1}{2}$, then
\begin{align}\label{kj2}
\sum_{k=j}^{n-1}q^{2k}
{k\brack j}_{q^2}^2
&\equiv\sum_{\alpha=0}^{\frac{n}{d}-1}{\alpha\choose\gamma}^2 \sum_{\beta=0}^{d-1}
q^{2\beta}{\beta\brack \delta}_{q^2}^2\notag\\
&=\sum_{\alpha=0}^{\frac{n}{d}-1}{\alpha\choose\gamma}^2 \sum_{r=0}^{d-\delta-1}
q^{2\delta+4r+4r\delta+2r^2} {-\delta-1\brack r}_{q^2}^2\notag\\
&\equiv \sum_{\alpha=0}^{\frac{n}{d}-1}{\alpha\choose\gamma}^2 q^{-2\delta-2\delta^2-2}{2d-2\delta-2\brack d-1-\delta}_{q^2}^2\pmod{\Phi_d(q)}.
\end{align}
It is obvious that the right-hand side of \eqref{kj2} divisible by $\Phi_d(q)$ for $\delta\leqslant \frac{d-3}{2}$, which means that
\begin{align*}\sum_{k=0}^{n-1}q^{2k}g_k(x;q^2)
&\equiv
\sum_{\gamma=0}^{\frac{n}{d}-1}\sum_{\alpha=0}^{\frac{n}{d}-1}{2\gamma\choose\gamma}{\alpha\choose\gamma}^2 x^{\gamma d+\frac{d-1}{2}}
q^{\frac{d-3}{2}} {d-1\brack \frac{d-1}{2}}_{q^2}^2\\
&\equiv \sum_{\gamma=0}^{\frac{n}{d}-1}\sum_{\alpha=0}^{\frac{n}{d}-1}{2\gamma\choose\gamma}{\alpha\choose\gamma}^2
 \frac{x^{\gamma d+\frac{d-1}{2}}}{q}  \pmod{\Phi_d(q)}.
\end{align*}
This proves \eqref{oddcong}.

Now assume that $d$ is an even divisor of $n$.
Similarly to \eqref{kj2}, we have
$$
\sum_{k=j}^{n-1}q^{k} {k\brack j}_{q}^2\equiv 0\pmod{\Phi_d(q)}\quad\text{for $j=\gamma d+\delta$ and $0\leqslant \delta\leqslant \frac{d}{2}-1$.}
$$
On the other hand, if $j=\gamma d+\delta$ with $\frac{d}{2}\leqslant \delta\leqslant d-1$, then by the $q$-Lucas theorem, we obtain
$$
{2j\brack j}_{q}={2\gamma d+2\delta\brack \gamma d+\delta}_{q}\equiv 0\pmod{\Phi_d(q)}.
$$
This proves \eqref{evencong1}.

Suppose that $d>2$ is even and $d\mid n$. Similarly to \eqref{oddqlucas} and \eqref{oddcv},
we get
$$
\sum_{k=j}^{n-1}q^k {k\brack j}_q {k+1\brack j+1}_q\equiv \sum_{\alpha=0}^{\frac{n}{d}-1}{\alpha\choose\gamma}^2
q^{-2\delta-\delta^2-2} {2d-2\delta-3\brack d-1-\delta}_{q}\equiv0\pmod{\Phi_d(q)}
$$
for $j=\gamma d+\delta$ and $0\leqslant \delta\leqslant \frac{d}{2}-2$.
On the other hand, if $j=\gamma d+\delta$ with $\frac{d}{2}\leqslant \delta\leqslant d-1$, then
$$
{2j\brack j}_q\equiv0\pmod{\Phi_d(q)}.
$$
while if $j=\gamma d+\frac{d}{2}-1$, then
$$
[j+1]_{q^2}\equiv\frac{1-q^d}{1-q^2}\equiv0\pmod{\Phi_d(q)}.
$$
This proves \eqref{evencong2}.  \qed

Recall that for $d>1$, we have
$$
\Phi_d(1)=\begin{cases} p,&\text{if }d=p^\alpha\text{ is a prime power},\\
1&\text{otherwise}.
\end{cases}
$$
Write $n=2^{r}n_1$, where $n_1$ is an odd integer. Then
$$
\prod_{\substack{d\mid n\\ d>1\text{ is odd}}}\Phi_d(1)=n_1,
\quad\text{and}\quad \prod_{\substack{d\mid n\\ d\text{ is even}}}\Phi_d(1)=2^{r}.
$$
Letting $q=1$ in \eqref{oddcong}--\eqref{evencong2}, we immediately get
\begin{align}
\sum_{k=0}^{n-1}(4k+3)g_k(x)\equiv 0\pmod{n_1}, \label{eq:final-1}
\end{align}
and
\begin{align}
2\sum_{k=0}^{n-1}kg_k(x)\equiv \sum_{k=0}^{n-1}g_k(x)\equiv 0\pmod{2^{r}}. \label{eq:final-2}
\end{align}
It is clear that \eqref{eq:first} follows from \eqref{eq:final-1} and \eqref{eq:final-2}. Therefore, Theorem \ref{thm:5}
may be deemed a $q$-analogue of \eqref{eq:first}.

\section{An open problem}
Numerical calculation suggests the following conjecture on congruences involving $S_n$.
\begin{conj}Let $n$ be a positive integer and $p$ a prime. Then
\begin{align}
\sum_{k=1}^{n}(-1)^k\frac{S_{k+2}+12S_{k+1}+16S_k}{k} &\equiv 0\pmod{n}, \label{eq:last-1}\\
\sum_{k=1}^{p}(-1)^k\frac{S_{k+2}+12S_{k+1}+16S_k}{k} &\equiv 2p(-1)^{\frac{p+1}{2}}\pmod{p^2}. \notag
\end{align}
\end{conj}

By \eqref{eq:multi-sum-3}--\eqref{eq:rewrite}, it is easy to see that if the $n=p$ case of \eqref{eq:last-1} is true, then we have
$$
\sum_{k=0}^{p-1}(8k^2+12k+5)g_k(-1) \equiv 0\pmod{p^2}.
$$
which is a special case of \eqref{eq:second-2} and \eqref{eq:third-2} conjectured by Z.-W. Sun.

\vskip 5mm \noindent{\bf Acknowledgment.} The authors would like to thank Professor Zhi-Wei Sun
for helpful comments.


\begin{thebibliography}{99}
\small \setlength{\itemsep}{-.8mm}

\bibitem{Andrews75}G.E. Andrews, A theorem on reciprocal polynomials with applications to
permutations and compositions, Amer. Math. Monthly 82 (1975), 830--833.

\bibitem{Andrews98}G.E. Andrews, The Theory of Partitions, Cambridge University Press, Cambridge, 1998.

\bibitem{AndrCB}
G.E. Andrews, The {F}riedman--{J}oichi--{S}tanton monotonicity
  conjecture at primes, Unusual Applications of Number Theory (M.~Nathanson,
  ed.), DIMACS Ser.\ Discrete Math. Theor. Comp. Sci., vol. 64, Amer. Math. Soc., Providence, R.I., 2004, pp.~9--15.

\bibitem{Des}J. D\'esarm\'enien, Un analogue des congruences de Kummer pour les $q$-nombres d'Euler,
European J. Combin. 3 (1982), 19--28.

\bibitem{Gessel}I.M. Gessel, Super ballot numbers, J. Symbolic Comput. 14 (1992), 179--194.

\bibitem{Guo}V.J.W. Guo, Proof of two divisibility properties of binomial coefficients conjectured by Z.-W. Sun,
Electron. J. Combin. 21(2) (2014), \#P2.54.

\bibitem{GK}V.J.W. Guo and C. Krattenthaler,
Some divisibility properties of binomial and $q$-binomial coefficients, J. Number Theory 135 (2014), 167--184.

\bibitem{GZ}V.J.W. Guo and J. Zeng, Some arithmetic properties of the $q$-Euler numbers and $q$-Sali\'e numbers,
European J. Combin. 27 (2006), 884--895.

\bibitem{Koepf}W. Koepf, Hypergeometric Summation, an Algorithmic Approach to Summation
and Special Function Identities, Friedr. Vieweg \& Sohn, Braunschweig, 1998.

\bibitem{MS}G.-S. Mao and Z.-W. Sun, Two congruences involving harmonic numbers with applications,
Int. J. Number Theory, in press.

\bibitem{Olive}G. Olive, Generalized powers, Amer. Math. Monthly 72 (1965), 619--627.

\bibitem{PWZ}M. Petkov\v{s}ek, H. S. Wilf and D. Zeilberger, $A=B$, A K Peters, Ltd., Wellesley, MA, 1996.

\bibitem{ReSWAA}
V. Reiner, D. Stanton, and D. White, The cyclic sieving phenomenon, J. Combin. Theory Ser. A 108 (2004), 17--50.

\bibitem{Riordan}J. Riordan, Combinatorial Identities, J. Wiley, New York, 1979.

\bibitem{Stanley89}R.P. Stanley, Log-concave and unimodal sequences in algebra, combinatorics, and geometry,
in: Graph Theory and Its Applications: East and West (Jinan, 1986),
Ann. New York Acad. Sci., 576, New York Acad. Sci., New York, 1989, pp.~500--535.

\bibitem{Stanley}R.P. Stanley, Enumerative Combinatorics,
  vol.~2, Cambridge University Press, Cambridge, 1999.

\bibitem{Sun00}Z.-W. Sun, On sums of Ap\'ery polynomials and related congruences, J. Number Theory
132 (2012), 2673--2699.

\bibitem{Sun13}Z.-W. Sun, Connections between $p=x^2+3y^2$ and Franel numbers, 
J. Number Theory 133 (2013), 2914--2928.

\bibitem{Sunnew}Z.-W. Sun, Conjectures and results on $x^2\mod p^2$ with $p^2=x^2+dy^2$, in: Number Theory and Related Area, Y. Ouyang, C. Xing, F. Xu
and P. Zhang, Eds., Adv. Lect. Math. 27, Higher Education Press and International Press, Beijing--Boston, 2013, pp. 149--197.

\bibitem{Sun}Z.-W. Sun, Congruences involving $g_n(x)=\sum_{k=0}^n{n\choose k}^2{2k\choose k}x^k$, Ramanujan J., doi:10.1007/s11139-015-9727-3.
\end{thebibliography}
\end{document}